\documentclass{amsart}

\usepackage[paperwidth=210mm,paperheight=297mm,hmargin={35mm,35mm},vmargin={35mm,35mm}]{geometry} %margens
\usepackage{amsthm,amssymb,latexsym,amsmath}
\usepackage{mathrsfs}
\usepackage{graphicx}
\usepackage{hyperref}
\input xy
\xyoption{all}
\usepackage[all]{xy}

\renewcommand{\dim}{\mathrm{dim}}

\newcommand{\Sing}{{\rm{Sing}}}
\newcommand{\Res}{{\rm{Res}}}
\newcommand{\codim}{{\rm{codim}}}

\newtheorem{nada}{Nada}[section]
\newtheorem{definition}[nada]{Definition}
\newtheorem{proposition}[nada]{Proposition}

\newtheorem{theorem}[nada]{Theorem}
\newtheorem*{thm*}{Theorem}
\newtheorem{lemma}[nada]{Lemma}

\newcommand{\bc}{\begin{center}}
\newcommand{\ec}{\end{center}}

\theoremstyle{plain}

\begin{document}

\title{A brief introduction on residue theory of holomorphic foliations}
\hyphenation{ho-mo-lo-gi-cal}
\hyphenation{fo-lia-tion}
\hyphenation{ge-ne-ra-li-za-ti-on}

%\dedicatory{\it Dedicated to Professor Jean-Paul Brasselet on the occasion of his 75th birthday}

\begin{abstract} This is a survey paper dealing with holomorphic foliations, with emphasis on residue theory and its applications. We start recalling the definition of holomorphic foliations as a subsheaf of the tangent sheaf of a manifold. The theory of Characteristic Classes of vector bundles is  approached from this perspective. We define Chern classes of holomorphic foliations using the Chern-Weil theory and we remark that the Baum-Bott residue is a great tool that help us to classify some foliations. We present throughout the survey several recent results and advances in residue theory. We finish by presenting some applications of residues to solve for example the Poincar\'e problem and the existence of minimal sets for foliations.

\end{abstract}

\author{Fernando Louren\c co}
\address{ Fernando Louren\c co \\ DMM - UFLA ,  Campus Universit\'ario, Lavras MG, Brazil, CEP 37200-000}
\email{fernando.lourenco@ufla.br}

\author{Fernando Reis}
\address{Fernando Reis \\ DMA - UFES,  Rodovia BR-101, Km 60, Bairro Litor\^aneo, S\~ao Mateus-ES, Brazil, CEP 29932-540}
\email{fernando.reis@ufes.br }

%\thanks{}
%\dedicatory{}
%\commby{ }
%\begin{abstract}
%\end{abstract}
%\begin{center}
\maketitle
%\end{center}
%\tableofcontents

\smallskip


\section{Introduction} % Please enter the title of your first section (only the first letter of the title should be capital).

The residue theory of holomorphic foliations red with the work of  P.Baum and R.Bott \cite{Article-2} in 1972. In their paper the authors have developed a class for foliations associated with its singular set using differential geometry based on the Bott vanishing theorem. However this class, called residue,  appears only as an element in the homology group. The question that arises is: \textit{how to  compute the residue?} We consider a holomorphic foliation $\mathcal{F}$ of dimension $q$ on a complex manifold $M$ of dimension $n$. If we take $\varphi$ a homogeneous symmetric polynomial of degree $d$, then to each compact connected component $Z \in \Sing(\mathcal{F})$ of the singular set of the foliation $\mathcal{F}$, there exists the homology class $\Res_{\varphi}(\mathcal{F};Z)$ into the group $H_{2(n-d)}(Z; \mathbb{C})$ such that under certain conditions on $M$, one has 

$$ \varphi(\mathcal{N}_{\mathcal{F}}) \frown [M] = \sum \Res_{\varphi}(\mathcal{F};Z), $$

\noindent where $\mathcal{N}_{\mathcal{F}}$ represents the normal sheaf of the foliation $\mathcal{F}
$.

This survey address the problem of computing this residue $\Res_{\varphi}(\mathcal{F};Z)$ in some cases. In one complex variable, we have the Cauchy's residue of a holomorphic function and the Cauchy integral formula which help us to compute it. On the other hand, in several complex variable, as a generalization of the Cauchy's residue, we have the Grothendieck residue associated with a meromorphic form.
Let $f = (f_{1}, \ldots , f_{n}) : U \subset \mathbb{C}^{n} \longrightarrow \mathbb{C}^{n}$ be a finite holomorphic map, such that $f(0) = 0,$ and $g$ be a holomorphic function on $U$. We define the Grothendieck residue by

$$\Res_{0}(g,f) = \Big(\dfrac{1}{2\pi \sqrt{-1}}\Big)^{n} \int_{\gamma} \dfrac{g(z) dz_{1}\wedge \ldots \wedge dz_{n}}{f_{1}\ldots f_{n}}$$

\noindent where $\gamma$ is a $n$-cycle with orientation prescribed by the $n$-form $d(\arg (f_{1})) \wedge \cdots \wedge d(\arg (f_{n}) )\geq 0.$ If we denote the merophorfic $n$-form $\dfrac{g(z) dz_{1}\wedge \cdots \wedge dz_{n} }{f_{1}\ldots f_{n}}$ by $\omega$ we may use the notation

\begin{equation}\label{eq 01(groth)}
\Res_{0}(g,f) = \Res_{0}\left[\begin{array}{cc} \omega \\

f_{1}, \ldots , f_{n}

\end{array}\right].
\end{equation}
We observe that for $n=1$, this residue is just the Cauchy's residue. Baum and Bott \cite{Article-2} also show how to compute the residue of a higher dimension foliation since certain assumptions are required in each irreducible component of the singular set of foliation. Take $Z$ an irreducible component of $\Sing(\mathcal{F})$ with $\dim Z = \dim (\mathcal{F}) - 1$ and other generic hypotheses, one has

$$\Res_{\varphi}(\mathcal{F}; Z) = \Res_{\varphi}(\mathcal{F}|_{B_{p}}; p)[Z],$$ 

\noindent where $\Res_{\varphi}(\mathcal{F}|_{B_{p}}; p)$ is a certain Grothendieck residue. 

Furthermore, in the same work, \cite{Article-2} Baum-Bott prove that the residue of a dimension one foliation at an isolated singular point can be also expressed by the Grothendieck residue in (\ref{eq 01(groth)}), where $f_{1}, \ldots , f_{n}$ are the components of the vector fields that locally induces $\mathcal{F}$.

In 1984, T. Suwa in \cite{Article-29} considers a foliation of complete intersection type and express a certain class of residues in terms of the Chern classes and the local Chern classes of the sheaf $\mathcal{E} xt_{\mathcal{O}_{M}}^{1}(\Omega_{\mathcal{F}}, \mathcal{O}_{M})$. As an application, in  (\cite{Article-29}, 3.8 Corollary, p.41) he gives a partial answer to {the Rationality Conjecture}, see (\cite{Article-2}, p.287). As another consequence, he shows how to compute residues at isolated singularities  in the case that the foliation has codimension one.

T. Suwa in \cite{book-5} has developed a residue theory to distributions, where the localization considered there arises from a rather primitive fact, i.e., the Chern forms of degree greater than the rank of the vector bundle vanish. Hence the involutivity has nothing to do with it. For this reason the localization can be defined by rank reason, of some characteristic classes and associated residues of the normal sheaf of the distribution. Also in  (\cite{book-5}, Proposition 4.4, p.15),  the author shows in particular that the residues to distribution coincides with the Baum-Bott residues to foliation, when the distribution is involutive (i.e., the distribution is integrable and hence define a foliation).

 A few years later, in 2005, based on the classical Camacho-Sad residue (or index) theorem, T. Suwa and F. Bracci in \cite{Article-6},  have developed a residue theory  for adequate singular pairs, which generalizes, in  some sense, the classical Camacho-Sad residue. The authors show a Bott vanishing theorem  by adapting  the Cech-de Rham theory and localization for adequate singular pair and prove that there exists the residue in this situation.

More recently, in 2015, F. Bracci and T. Suwa in \cite{Article-7} provide another way to compute residues. The authors consider a deformation of a complex manifold $M$, denoted by $\tilde{M} = \{ M_{t} \}_{t \in \hat{P}}$, where the parameter space $\hat{P}$ is a $\mathcal{C}^{\infty}$ manifold and a deformation of a holomorphic foliation $\mathcal{F}$ of $M$, denoted by $\mathcal{\tilde{F}} = \{ \mathcal{F}_{t} \}$. For each parameter $t \in \hat{P}$ they assume that the singular set $S_{t}$ is compact. Let $\varphi$ be a homogeneous symmetric polynomial of degree $d$ and $\Res_{\varphi}(\mathcal{F}_{t}; S_{t}^{\lambda})$ the Baum-Bott residue of  $\mathcal{F}_{t}$ at the compact connected set $S_{t}^{\lambda}$. It is proved that the Baum-Bott residue continuously  varies on this smooth deformation, i.e,

$$\lim_{t \rightarrow t_{0}} \sum_{\lambda}\Res_{\varphi}(\mathcal{F}_{t}; S_{t}^{\lambda}) = \Res_{\varphi}(\mathcal{F}_{t_{0}}; S_{t_{0}}^{\lambda}).$$ 
Subsequentely, the authors consider $\mathcal{F}$ as a germ at $ 0 \in \mathbb{C}^{n}$ of a simple almost Liouvillian foliation of codimension one and $V$ a divisor of poles. Then it is shown that the residue of $\mathcal{F}$ at $Z$, which is an irreducible component of singular set of $\mathcal{F}$ of codimension 2, can be written as a sum (over the irreducible components of $V$ that contains $Z$) in terms of Lehmann-Suwa residues \cite{Article-23}. This represents an effective way to compute residues:

$$BB(\mathcal{F}; Z) = \sum_{j=1}^{k} \Res(\gamma_{0},V_{j})Var(\mathcal{F},V_{j};Z).$$

In the paper \cite{Article-31} the author proves a more  slight generalization of the Bott residue theorem to holomorphic foliations of dimension one. The proof is based on a localization formula of Duistermaat and Heckman type, which has been first discussed in \cite{Article-5}.
 
Recently in 2016, see \cite{Article-14.1},  Corr\^ea, Pen\~a and Soares have studied several residue formulas for vector fields on compact complex orbifolds with isolated singularities,  which is a special type of a singular variety.

It is worth remark that there are other types of residues and invariants associated with a foliations and distributions. For instance, residues of logarithmic vector fields in \cite{Article-14.3}, Camacho-Sad index in \cite{Article-11},  GSV-index of foliations and Pfaff systems in \cite{Article-14.2,Article-20.1}. Another interesting topic that has been studied recently is the residue theory for  flags of foliations and distributions. In the works \cite{Article-9, Article-20}, it was developed a general theory of residues to $2$-flags. There are many topics which are closely related to $2$-flags
and naturally appear  in the theory of foliation. For example,
	a conjecture due to Marco Brunella (see \cite{Article-12}, p.443) which says that a two-dimensional holomorphic foliation $\mathcal{F}$ on $\mathbb{P}^{3}$ either admits an invariant algebraic surface or it compose a flag of holomorphic foliations. 

 Finally, we finish the survey with some applications. We discuss how the residue theory has been used to investigate two classical open problems: the Poincar\'e problem and the existence or not of minimal sets for foliations.

It is important to note that throughout this paper we maintain all the notations and symbols exactly as they are in the original works, respecting the choices of the referenced authors.

\section{Chern-Weil Theory of Characteristic classes of holomorphic foliations} % Please enter the title of your second section.

\subsection{Holomorphic foliations}

Let us begin by recalling the basic definition of singular holomorphic foliations and distributions. Let $M$ be a complex manifold of dimension $n$ and denote by $\Theta_{M}$ and $\Omega_{M}$ respectively, the sheaves of germs of holomorphic vector fields and holomorphic $1$-forms on $M$. There are two definitions of singular foliations that turn out to be equivalent as long as we consider only reduced foliations. For this section about foliations theory we refer to  \cite{Article-2, Article-9, Article-24, book-4, book-5}.

A singular holomorphic distribution of dimension $q$ on $M$ is a coherent subsheaf $\mathcal{F}$ of $\Theta_{M}$ of rank $q$. Moreover, if $\mathcal{F}$ satisfies the following integrability condition

$$ [\mathcal{F}_{x} , \mathcal{F}_{x}] \subset \mathcal{F}_{x} \ \ \ \mbox{for} \ \  \mbox{all} \ \ x \in M,$$

\noindent we say that $\mathcal{F}$ is a holomorphic foliation. The normal sheaf of $\mathcal{F}$ is defined as the quotient sheaf $\mathcal{N}_{\mathcal{F}} := \Theta_{M}/ \mathcal{F}$, such that it is torsion free (it means that $\mathcal{F}$ is
saturated).
With this definition we have the following exact sequence 
$$  0  \longrightarrow \mathcal{F}  \longrightarrow  \Theta_{M}  \longrightarrow \mathcal{N}_{\mathcal{F}} \longrightarrow 0.
$$
We define the singular set of the distribution $\mathcal{F}$ by  
\begin{center}
    $\Sing(\mathcal{F}) := \Sing(\mathcal{N}_{\mathcal{F}} ) = \{p \in M ; \mathcal{N}_{\mathcal{F},p} \  \mbox{is} \ \ \mbox{not} \  \mbox{locally}   \ \mbox{free} \}.$
\end{center}
We assume that $\codim(\Sing\mathcal{(F)}) \geq 2.$

For the second one definition, a singular distribution $\mathcal{G}$ can be defined, as a dual way by means of differential forms, i. e., as a coherent subsheaf of $\Omega_{M}$. If $\mathcal{G}$ satisfies the integrability condition, i.e.,

$$d \mathcal{G}_{x} \subset (\Omega_{M} \wedge \mathcal{G})_{x} \ \ \mbox{for} \ \  \mbox{all}  \ x \in M \setminus \Sing(\mathcal{G}),$$

\noindent we say that $\mathcal{G}$ is a foliation, where  $\Sing(\mathcal{G}) := \Sing(\Omega_{M}
/ \mathcal{G}).$

The two definitions of foliations are equivalents and related by taking the annihilator of each other. If $\mathcal{F}$ is a foliation on $M$ of dimension $q$, its annihilator is defined by

$$ \mathcal{F}^{a} = \{ v \in \Omega_{M} ; < v, \omega > = 0 \ \ \mbox{for} \ \  \mbox{all} \ \ \omega \in \mathcal{F} \}.$$

We say that $\mathcal{F}$ is reduced, if for any open set $U$ in $M$,

$$\Gamma(U, \Theta_{M}) \cap \Gamma(U\setminus \Sing(\mathcal{F}), \mathcal{F}) = \Gamma(U, \mathcal{F}).  $$

If we consider only reduced foliation, then $\mathcal{G} = \mathcal{F}^{a}$ and the converse is also true (see \cite{book-4}).

To conclude this subsection we present the follow definition.

\begin{definition} Let $V$ be an analytic subspace of a complex manifold $X$. We say that $V$ is invariant by a foliation $\mathcal{F}$ if

$$T\mathcal{F}|_{V} \subset (\Omega_{V}^{1})^{\ast}.$$

\noindent In particular cases,

$\bullet$ if $V$ is a hypersurface we say that $\mathcal{F}$ is \textit{logarithmic \ along} $\mathcal{F}$;

$\bullet$ if $V$ is a reduced complete intersection of dimension $n-k$, defined by intersection of $k$ hypersurfaces we say that $\mathcal{F}$ is \textit{multlogarithmic \ along} $\mathcal{F}$.
\end{definition}
\subsubsection{Flag of holomorphic foliations}

In this subsection, we define flag of foliations and show its main properties. For more details we refer to \cite{Article-9, Article-12.22, Article-20, Article-24}. 

A flag of singular holomorphic foliations on a complex manifold $M$ of dimension $n$, can be define by a finite sequence of $k$ foliations $\mathcal{F} = (\mathcal{F}_{1},\ldots , \mathcal{F}_{k})$ such that, outside the singular sets, each foliation $\mathcal{F}_{i+1}$ is a  subfoliation of $\mathcal{F}_{i}$ and we denote $\mathcal{F}_{i} \subset \mathcal{F}_{i+1}$, for each $i = 1, \ldots, k-1$. In a more formal manner, and for $k=2$, one has the following. 

%\vspace{-1.5 cm}

\begin{definition}\label{def.2.1} Let $\mathcal{F}_{1},\mathcal{F}_{2}$ be two holomorphic foliations on $M$ of dimensions $q = (q_{1},q_{2})$. We say that $\mathcal{F} := (\mathcal{F}_{1},\mathcal{F}_{2})$ is a $2$-flag of holomorphic foliations if $\mathcal{F}_{1}$ is a coherent
sub $\mathcal{O}_{M}$-module of $\mathcal{F}_{2}.$
\end{definition}

%\vspace{-0.18 cm}
We note that, for $ x \in M \setminus \cup_{i = 1}^{2} \Sing (\mathcal{F}_{i})$ the inclusion relation $T_{x}\mathcal{F}_{1} \subset T_{x}\mathcal{F}_{2}$ holds, namely, the leaves of $\mathcal{F}_{1}$ are contained in leaves of $\mathcal{F}_{2}$. Here $T\mathcal{F}_{i}$ represents the tangent sheaf of the foliation $\mathcal{F}_{i}$, but throughout the text we will abuse notation and denote it simply by $\mathcal{F}_{i}$. Now, we would like to highlight a diagram of short exact sequences of sheaves, called "turtle diagram".

$$ \xymatrix{ 0 \ar[rd]  &  & 0 \ar[ld] &  & 0  \\
   & \mathcal{F}_{1} \ar[rd] \ar[dd] &  & \mathcal{N}_{2} \ar[lu] \ar[ru] &  \\
   &  & \Theta_{M} \ar[ru] \ar[rd]  &  &  \\
   & \mathcal{F}_{2} \ar[ru]   \ar[rd] &  & \mathcal{N}_{1} \ar[uu]  \ar[rd] &  \\
 0 \ar[ru] &  & \mathcal{N}_{12} \ar[ru] \ar[rd]  &  & 0 \\
  & 0 \ar[ru] &   &  0 &  }$$
Let us give some definitions and notations. We write $\mathcal{N}_{j}$ as the quotient sheaf  $\mathcal{N}_{\mathcal{F}_{j}}$. Also, we denote by $\Sing(\mathcal{F})$ the singular set of the flag $\mathcal{F}$ which is the analytic set $\Sing(\mathcal{F}_{1})\cup \Sing(\mathcal{F}_{2})$, and $\mathcal{N}_{\mathcal{F}} := \mathcal{N}_{12} \oplus
\mathcal{N}_{2}$ is the normal sheaf of the flag, where
$\mathcal{N}_{12}$ is the relative quotient sheaf $ \mathcal{F}_{2} / \mathcal{F}_{1}$.

\subsection{Characteristic classes via Chern-Weil theory}

In this section, we review the basic tools of the Chern-Weil theory for working with residue and characteristic classes to vector bundles and sheaves. The residue theory of foliations was first introduced by Baum and Bott using differential geometry, in a series of papers (\cite{Article-2, Article-3, Article-4}). Later Lehman and Suwa in the decades of 1980 and 1990 present a new approach of residue theory using the Chern-Weil theory (see \cite{Article-23}).

\begin{definition} A connection for a complex vector bundle $E$ of rank $r$ on $M$ is a $\mathbb{C}$-linear map
$$ \nabla : A^{0} (M,E) \longrightarrow A^{1} (M, E) $$

\noindent that satisfies

$$ \nabla (f.s) = df \otimes s + f. \nabla (s) \ \ \mbox{for} \ \ f \in A^{0}(M) \ \ \mbox{and} \ \ s \in A^{0}(M,E).$$

\end{definition}
If $\nabla$ is a connection for $E$, then it induces a $\mathbb{C}$-linear map
$$ \nabla := \nabla^{2} : A^{1}(M,E) \longrightarrow A^{2}(M,E) $$

\noindent satisfying
$$ \nabla (\omega \otimes s) = d\omega \otimes s - \omega \wedge \nabla (s), \ \ \omega \in A^{1}(M), \ \ s \in A^{0} (M,E).$$

\noindent We define the composition map $K := \nabla \circ \nabla$ from $A^{0} (M,E)$ to $A^{2} (M,E)$ as the curvature of the connection $\nabla$. If $ s = (s_{1}, \ldots , s_{n})$ is a frame of $E$ on an open set $U$ we have $\theta = (\theta_{ij})$ the connection matrix (where $\theta_{ij}$ are 1-forms on $U$ ) of $E$ with respect to frame $s$. In the same way, we can get $K = (k_{ij})$ the curvature matrix of $E$ with respect to $s$. If we consider $\sigma_{i}, i = 1,\ldots,r$ the $i$-th elementary symmetric functions in the eigenvalues of the matrix $K$
$$ \det(It + K) = 1 + \sigma_{1}(K)t + \sigma_{2}(K)t^{2} + \cdots + \sigma_{r}(K)t^{r}, $$

\noindent we may define a $2i$-form of Chern $c_{i}$ on $U$ by
$$ c_{i}(K) := \sigma_{i} (\frac{\sqrt{-1}}{2 \pi}K). $$

In general, if $\varphi$ is a homogeneous symmetric polynomial in $r$ variables of degree $d$, we may write $\varphi = \tilde{P}(c_{1},\ldots,c_{r})$ for some polynomial $\tilde{P}$. Then we can define
$$ \varphi(K) := \tilde{P}(c_{1}(K),\ldots,c_{r}(K))$$
\noindent which is a closed form on $M$. Therefore, we have a cohomology class of $E$ on $M$, $\varphi (E) := \varphi(K) \in H^{2d}(M; \mathbb{C})$.

Let $\mathcal{G}$ be a sheaf on $M$, $S$ a compact connected set of $M$ and $U$ a relatively compact open neighborhood of $S$ in $M$. We may consider $\mathcal{U} = \{U_{0}, U_{1} \}$ a covering of $U$, where $U_{1} = U$ and $U_{0} = U \setminus S$ and since there exists \cite{Article-1} the following resolution of $\mathcal{G}$

$$0 \longrightarrow \mathcal{A}_{U}(E_{r}) \longrightarrow \cdots \longrightarrow \mathcal{A}_{U}(E_{0}) \longrightarrow \mathcal{A}_{U} \otimes_{\mathcal{O}_{U}}\mathcal{G} \longrightarrow 0,$$

\noindent we can define the characteristic class $\varphi(\mathcal{G})$ on $U$ using the \textit{virtual bundle} $\xi := \sum_{i=0}^{r} (-1)^{i}E_{i}$, i. e.,
$\varphi(\mathcal{G}) :=\varphi(\xi).$

Given $\mathcal{F}$ a holomorphic foliation on $M$ and $\varphi$ a homogeneous symmetric polynomial of degree $d$, one has the short exact sequence

$$ 0 \longrightarrow \mathcal{F} \longrightarrow \Theta_{M} \longrightarrow \mathcal{N}_{\mathcal{F}} \longrightarrow 0.$$

Then $\varphi(\mathcal{N}_{\mathcal{F}})$ denotes the characteristic class of $\mathcal{F}$ and it is an element in cohomology group $H^{2d}(M; \mathbb{C}).$ We denote by $P$ the Poincar\'e homomorphism (or isomorphism if $M$ is nonsingular) from $H^{2d}(M;\mathbb{C})$ to $H_{2(n-d)}(M;\mathbb{C})$ and by $A$ the Alexander homomorphism (or isomorphism if $S$ is nonsingular) $A : H^{2d}(M, M \backslash S; \mathbb{C}) \longrightarrow H_{2(n-d)}(S; \mathbb{C})$. We have the
following commuting diagram:

$$ \xymatrix{ H^{2d}(M, M \backslash S; \mathbb{C}) \ar[r]  \ar[d]_{A} & H^{2d}(M ; \mathbb{C})  \\ H_{2(n-d)}(S; \mathbb{C}) \ar[r]^{i^{\ast}} & H_{2(n-d)}(M; \mathbb{C}) \ar[u]_{P}}$$

\noindent where this map $H^{2d}(M, M \backslash S; \mathbb{C}) \longrightarrow H^{2d}(M ; \mathbb{C}) $ represents a lift that can be interpreted, in terms of foliation theory, by the Bott vanishing theorem (see \cite{book-4}, Theorem 9.11). Thus we have the residue of foliation $\mathcal{F}$, denoted by $\Res_{\varphi}(\mathcal{F}, \mathcal{N}_{\mathcal{F}}; S)$ in $H_{2(n-d)}(S; \mathbb{C})$ as the image of $\varphi_{S}(\mathcal{N}_{\mathcal{F}}; \mathbb{C}) \in H^{2d}(M, M \backslash S; \mathbb{C})$ by the Alexander duality, which is independent of all the choices involved.

In general it is not possible to compute such residue directly from the above definition. It is then important when one can compute such element using tools like differential geometry, foliation theory, complex analysis or singularities theory. The goal of this survey is to present some results in this direction.

\subsection{Some results about residues of holomorphic foliations}

The residue theory of holomorphic foliations was developed by several authors in the last years. We can cite, for instance, Baum and Bott in \cite{Article-2, Article-4}, Brasselet and Suwa in \cite{Article-8} and Bracci and Suwa in \cite{Article-6} and \cite{Article-7}. We also refer to the book \cite{book-0} that is dedicated to study of indices for holomorphic foliation of dimension one with isolated singularity, in the cases where the underlying space is either a smooth variety or a singular variety. The authors defines several notions of index such that: Poincar\'e-Hopt index, Schwartz index, GSV index, Virtual index, Homological index and others.

We begin with the classical Grothendieck residue, see (\cite{book-2, Article-28, book-4}) for more details. Let us take $\omega$ a germ of holomorphic $n$-form at $0$, a neighborhood $U$ of $0$ in $\mathbb{C}^{n}$, and $a_{1}, \ldots , a_{n}$ germs of holomorphic functions such that $V(a_{1}, \ldots , a_{n}) = \{0\}$. The Grothendieck residue of $\omega$ at 0 is defined by

$$
\Res_{0}\left[\begin{array}{cc} \omega \\

a_{1}, \ldots , a_{n}

\end{array}\right] =  \dfrac{1}{(2\pi \sqrt{-1})^{n}} \int_{\gamma} \dfrac{\omega}{a_{1}\cdots a_{n}},
$$
\\
\noindent where $\gamma$ is a $n$-cycle in $U$ defined by
$$\gamma = \{ z \in U ; |a_{1}(z)| = \cdots = |a_{n}(z)|= \epsilon  \}$$

\noindent and oriented by $d(\arg( a_{1})) \wedge \cdots \wedge d(\arg( a_{n})) \geq 0.$ We remark that the above residue is the usual Cauchy residue at $0$, of the meromorphic 1-form $\omega /a_{1}$, when $n=1$.

Next, we give the first important result about residues in theory of foliations. This theorem relates the Baum-Bott residue of certain foliations, to the Grothendieck residue. For this, let $\mathcal{F}$ be a holomorphic foliation of dimension one in a complex compact manifold $M$ of dimension $n$. We assume that $\mathcal{F}$ has only isolated singularities. Let $\varphi$ be a homogeneous symmetric polynomial of degree $n$ and $p \in \Sing(\mathcal{F})$ an isolated point. Then

\begin{theorem}(\cite{Article-2}, Theorem 1)

$$\Res_{\varphi}(\mathcal{F}; p) = \Res_{p}\left[\begin{array}{cc} \varphi(J X) \\

X_{1}, \ldots , X_{n}

\end{array}\right],$$ 

\noindent where $X = (X_{1}, \ldots , X_{n})$ is a germ of a holomorphic
vector field at $p$, local representative of $\mathcal{F}$.

\end{theorem}

In \cite{Article-29}, (3.12) Proposition, T. Suwa considers $\mathcal{F}$ a holomorphic distribution (not necessarily involutive) of codimension one and taking $0$ (it can be another point $p$) as an isolated singularity of this distribution, he shows how to compute the residues.

\begin{theorem} Let $U$ be a polydisk about the origin $0$ in $\mathbb{C}^{n}$ and let $\mathcal{F} = <\omega>$ be a codimension one holomorphic foliation on $U$ with an isolated singularity at 0. Then we have

$$ \Res_{c_{n}}(\mathcal{F};0) = (-1)^{n}(n-1)! \dim_{\mathbb{C}}Ext^{1}_{\mathcal{O}}(\Omega_{\mathcal{F}}, \mathcal{O}) \ \ \ in \ \ \ H_{0}(0; \mathbb{Q}) = \mathbb{Q}.$$

\end{theorem}

The reader interested in residues for distributions can also refer the survey \cite{book-5}.

With the same aim, F. Bracci and T. Suwa in \cite{Article-7} consider smooth deformations of holomorphic foliations and verify that it provides an effective way to compute residues.

\begin{theorem} Let $(\tilde{M}, \hat{P}, \pi)$ be a deformation of manifolds and $\tilde{\mathcal{F}}$ a deformation of foliations on $\tilde{M}$ of rank $p$. Suppose that $\mathcal{N}_{\tilde{\mathcal{F}}}$ admits a $\mathcal{C}^{\infty}$ locally free resolution. Let $S^{'}(\tilde{\mathcal{F}}) \subset S(\tilde{\mathcal{F}})$ be a connect component of the singular set of $\tilde{\mathcal{F}}$ and let $S_{t} :=M_{t} \cap S^{'}(\tilde{\mathcal{F}})$. Assume that for all $t \in \hat{P}$ the set $S_{t}$ is compact and $S_{t} \neq M_{t}$. Let $\varphi$ be a homogeneous symmetric polynomial of degree $d > n-p$. Under this assumptions, the Baum-Bott residue $BB_{\varphi}(\mathcal{F}_{t};S_{t})$ is continuous in $t \in \hat{P}$. Namely, for any $\mathcal{C}^{\infty} \ (2n-2d)$-form $\tilde{\tau}$ on $\tilde{M}$ such that $i_{t}^{\ast}(\tilde{\tau})$ is closed for all $t \in \hat{P}$,

$$\lim_{t \rightarrow t_{0}} BB_{\varphi} (\mathcal{F}_{t};S_{t})(i_{t}^{\ast}(\tilde{\tau})) = BB_{\varphi} (\mathcal{F}_{t_{0}};S_{t_{0}})(i_{t_{0}}^{\ast}(\tilde{\tau})).$$

\end{theorem}

Also for higher dimension foliations it is possible to relate the Baum-Bott residue with the Grothendieck residue. Vishik, in \cite{Article-30}, found this relation with the hypothesis that the foliation has locally free tangent sheaf. In \cite{Article-2} Baum and Bott, before and independent of Vishik, have proved a similar result using a generic assumption in the singular set of foliation.

Let us consider $\mathcal{F}$ be a holomorphic foliation of codimension $k$ on a complex manifold $M$ and $\varphi$ a homogeneous symmetric polynomials of degree $k+1$. Note that $\deg \varphi > n - \dim (\mathcal{F})$, which is condition to Bott vanishing theorem. We consider that the singular set of $\mathcal{F}$ has pure expected codimension, i.e., $\dim (\Sing(\mathcal{F})) = k+1$.
In this case,  the notation $\Sing_{k+1}(\mathcal{F})$ is usually used to denote the union of irreducible components of the singular set of pure codimension $k+1$ of $\mathcal{F}$.

If $ Z \subset \Sing_{k+1}(\mathcal{F})$ is a pure codimension and irreducible component, we consider $B_{p}$ a $(k+1)$-dimension ball centered at $p$ sufficiently small and transversal to $Z$ at $p$. We remark that the restricted foliation $\mathcal{F}|_{B_{p}}$ is a foliation of dimension one with isolated singularity at $p$. In \cite{Article-14} (Theorem 1.2), M. Corr\^ea and F. Louren\c co relate the Baum-Bott residue of $\mathcal{F}$ in $Z$ with the Grothendick residue of $\mathcal{F}|_{B_{p}}$ at $p$.

\begin{theorem} Let $\mathcal{F}$ be a singular holomorphic foliation of codimension $k$ on a compact complex manifold $M$ such that $cod(\Sing(\mathcal{F} )) \geq k + 1.$ Then, 

$$\Res_{\varphi}(\mathcal{F}; Z) = \Res_{\varphi}(\mathcal{F}|_{B_{p}}; p)[Z],$$

\noindent where $\Res_{\varphi}(\mathcal{F}|_{B_{p}}; p)$ represents the Grothendieck residue at $p$ of the one dimensional foliation $\mathcal{F}|_{B_{p}}$ on a $(k + 1)$-dimensional transversal ball $B_{p}.$

\end{theorem}

The next result is due to Fernandez-Perez and Tamara in \cite{Article-19} (Theorem 6.2), and it provides another effective way to computing Baum-Bott residues of codimension one holomorphic foliations. First, we need a definition. Consider the germ of holomorphic foliation $\mathcal{F}$ at $0\in \mathbb{C}^n$.

\begin{definition} We say that $\mathcal{F}$ is an \textit{almost Liouvillian foliation} at $0 \in \mathbb{C}^{n}$ if there exists a germ of closed meromorphic $1$-form $\gamma_{0}$ and a germ of holomorphic $1$-form $\gamma_{1}$ at $0 \in \mathbb{C}^{n}$ such that

$$d \omega = (\gamma_{0} + \gamma_{1}) \wedge \omega.$$

\noindent We say that $\mathcal{F}$ is a \textit{simple almost Liouvillian foliation} at $0 \in \mathbb{C}^{n}$ if we can choose $\gamma_{0}$ having only first-order poles.

\end{definition}

\begin{theorem} Let $\mathcal{F}$ be a germ at $0 \in \mathbb{C}^{n} , n \geq 3$, of a simple almost Liouvillian foliation defined by $\omega \in \Omega^{1}(\mathbb{C}^{n},0)$ such that

$$ d\omega = (\gamma_{0}+ \gamma_{1}) \wedge \omega.$$

\noindent Let $V$ be the divisor of poles of $\gamma = \gamma_{1} + \gamma_{1}$ and $V_{1}, \ldots , V_{l}$ the irreducible components of $V$. Let $Z$ be an irreducible component of $\Sing_{2}(\mathcal{F})$. Then

$$ BB(\mathcal{F}; Z) = \sum_{i=1}^{k} \Res(\gamma_{0}, V_{j}) Var(\mathcal{F}, V_{j}; Z), $$

\noindent where $V_{1}, \ldots , V_{k}$ are the irreducible components that contains $Z$ and \\ $Var(\mathcal{F}, V_{j}; Z)$ represents the Varational index defined by Khanedani and Suwa in \cite{Article-21}.

\end{theorem}

In \cite{Article-8, Article-27} the notion of Nash residue of foliations was developed, which immediately implies a partial answer to the rationality conjecture of Baum and Bott (see \cite{Article-2}, p.287). Let $M$ be a complex manifold with dimension $n$, and let $\mathcal{F}$ be a singular holomorphic foliation with dimension $q$ on $M$. Let us consider for each point $x$ in $M$ the following set

\begin{equation}\label{eq1}
F(x) := \{ v(x); v \in \mathcal{F}_{x} \},
\end{equation}

\noindent where $\mathcal{F}_{x}$ denotes the stalk of $\mathcal{F}$ at $x$. We observe that $F(x)$ is a subspace of tangent space $T_{x}M$ of dimension $q$ if, and only if, $x \in M \setminus \Sing(\mathcal{F})$. In general, $\dim( F(x)) \leq q$. In the following, $G(q,n)$ is the Grassmannian bundle of $q$-planes in $TM$.

Using the expression $(\ref{eq1})$ we can define a section of $G(q,n)$ outside of singular set of $\mathcal{F}$, as follows

$$s: M \setminus \Sing(\mathcal{F}) \longrightarrow G(q,n)$$

\noindent gives by $s(x) := F(x).$

We define $M^{\nu} := \overline{Im (s)}$ in $G(q,n)$ and call it the Nash modification of $M$ with respect to foliation $\mathcal{F}$. We consider $\Sing(\mathcal{F})^{\nu} := \pi^{-1}\Sing(\mathcal{F})$ where $\pi$ is the restriction map to $M^{\nu}$ of the projection of the bundle $G(q,n)$ that is a birational map

$$\pi : M^{\nu} \longrightarrow M.$$

\noindent Moreover, it is biholomorphic from $M^{\nu} \setminus \Sing(\mathcal{F})^{\nu}$ to $M \setminus \Sing(\mathcal{F})$. In some case, we can assume $M^{\nu}$ as a smooth manifold (see \cite{Article-27}). We denote by $\tilde{T}^{\nu}$ and $\tilde{N}^{\nu}$, respectively, the tautological bundle and the tautological quotient bundle on $G(q,n)$. So, one has a short exact sequence

$$ 0 \longrightarrow T^{\nu} \longrightarrow \pi^{\ast}TM \longrightarrow N^{\nu} \longrightarrow 0, $$

\noindent where $T^{\nu}$ and $N^{\nu}$ are essentially the restrictions to $M^{\nu}$.

It is possible to show that the characteristic class $\varphi(N^{\nu})$, for a homogeneous symmetrical polynomial $\varphi$ of degree $d > n- \dim (\mathcal{F})$, is localized at $\Sing(\mathcal{F})^{\nu}$ give us the following residues

$$\Res_{\varphi}(N^{\nu}, \mathcal{F}; \Sing(\mathcal{F})^{\nu} ) \in H_{2(n-d)}(\Sing(\mathcal{F})^{\nu}; \mathbb{C})$$

\noindent and we call it the Nash residue of $\mathcal{F}$ with respect to $\varphi$ at $\Sing(\mathcal{F})^{\nu}$.

In 1989, Sert\"oz  in \cite{Article-27} (Theorem IV.4, p.238), has showed that the difference between the Baum-Bott residue and the Nash residue is an integer number with the assumption that $M^{\nu}$ is non singular.

\begin{theorem} Let $S$ be a connect component of singular set of $\mathcal{F}$ and $\varphi$ a homogeneous symmetrical polynomial of degree $d > n- \dim (\mathcal{F})$ then

$$\Res_{\varphi}(N_{\mathcal{F}}, \mathcal{F}; S) = \Res_{\varphi}(N^{\nu}, \mathcal{F}; S^{\nu}) +k,$$

\noindent where $k$ is a homology cycle in $S$ and it is compute by a Grassmaann graph construction.

\end{theorem}

In 2000, Brasselet and Suwa in \cite{Article-8} (Theorem 4.1, p. 44), have given a similar result of Sert\"oz droping the hypothesis that $M^{\nu}$ is smooth.

\begin{theorem} Let $\varphi$ be a homogeneous symmetric polynomial of degree $d > n - \dim (\mathcal{F})$. If $\varphi$ is with integral coefficients, then the difference

$$\Res_{\varphi}(N_{\mathcal{F}}, \mathcal{F}; S) - \Res_{\varphi}(N^{\nu}, \mathcal{F}; S^{\nu})$$

\noindent is in the image of the canonical homomorphism $$H_{2(n-p)}(S; \mathbb{Z}) \longrightarrow H_{2(n-p)}(S; \mathbb{C}).$$

\end{theorem}

In the following result F. Bracci and T. Suwa have developed the residue theory for foliations of adequate singular pairs (see \cite{Article-6}). In short, we consider $M$ a complex manifold of dimension $m$ and let $P \subset M$ be a complex submanifold of dimension $r$, then we have a short exact sequence

$$ 0 \longrightarrow TP \longrightarrow TM|_{P} \longrightarrow N_{P,M} \longrightarrow 0,$$

\noindent where $N_{P,M}$ denotes the normal bundle of $P$ in $M$. We pick $X$ another submanifold of $M$ of dimension $n$ which intersects $P$ along a submanifold $Y \subset M$ of dimension $n+r-m$ and such intersection is everywhere transversal. We define $(X, Y)$ as adequate singular pair in $M$ if $r = m+l-n$ and the data satisfy the following \\

1) $Y = X \cap P;$

2) $\dim (\Sing(X) \cap P) < l$;

3) $X_{reg}$ intersects $P$ generically transversely. \\

\begin{theorem}(\cite{Article-6}, Theorem 2.1, p.7) Let $(X,Y)$ be an adequate singular pair in $M$ and let  $\mathcal{F}$ be a holomorphic foliation in $X$ of dimension $d \leq l$ which leaves $Y$ invariant. Let $\Sigma = (\Sing(\mathcal{F}) \cup \Sing(Y)) \cap Y$ and assume that $\dim (\Sigma) < l$. Let $\Sigma = \sum_{\gamma} \Sigma_{\gamma}$ be the decomposition into connected components and let $ i_{\gamma} : \Sigma_{\gamma} \hookrightarrow Y$ denotes the inclusion. Let $\varphi$ be a symmetric homogeneous polynomial of degree $t>l-d$. Then

\begin{enumerate}

  \item[i)] For each compact connected component $\Sigma_{\gamma}$ there exists a class \\ $\Res_{\varphi}(\mathcal{F}, Y; \Sigma_{\gamma}) \in H_{2l-2t}(\Sigma_{\gamma}; \mathbb{C})$ called "residue", which depends only on the local behavior of $\mathcal{F}$ near $\Sigma_{\gamma}$; 
  
  \item[ii)] If $Y$ is compact we have
  
  $$ \sum_{\gamma} (i_{\gamma})_{\ast} \Res_{\varphi}(\mathcal{F}, Y; \Sigma_{\gamma}) = \varphi(N_{P,M}) \cap [Y] \ \ in \ \ H_{2l-2t}(Y; \mathbb{C}).$$
\end{enumerate}

\end{theorem}

We note that this "new concept" of residue of foliations can be understood as a generalization of the classical Camacho-Sad residue theorem  (see \cite{Article-11}).

Corr\^ea at al in \cite{Article-14.1} showed the follows residue formula for orbifolds. Let $X$ be a complex orbifold of dimension $n$ and $L$ be a line $V$-bundle over $X$ and considering some Chern classes of the bundle $TX-L^{\vee}$, moreover, for each point $p$ which vanishing $\xi$, let

$$\pi_{p} : (\tilde{U}, \tilde{p}) \rightarrow (U,p),$$

\noindent be a smoothing covering of $X$ at $p$ and the notation $\tilde{\xi} = \pi_{p}^{\ast} \xi$, one has

\begin{theorem}(\cite{Article-14.1}, Theorem 3.1, p.2897) Let $X$ be a compact orbifold of dimension $n$ with only isolated singularities, let $L$ be a locally $V$-free sheaf of rank $1$ over $X$ and $L$ the associated line $V$-bundle. Suppose $\xi$ is a holomorphic section of $TX\otimes L$ with isolated zeros. If $P$ is an invariant
polynomial of degree $n$, then

$$\int_{X} P(TX- L^{\vee}) = \sum_{p | \xi(p) = 0} \dfrac{1}{\#G_{p}} \Res_{\tilde{p}}\Big[P(J \tilde{\xi}) \dfrac{d\tilde{z}_{1}\wedge \ldots \wedge d\tilde{z}_{n}}{\tilde{\xi}_{1}\ldots \tilde{\xi}_{n}}\Big],$$

\noindent where $G_{p} \subset Gl(n, \mathbb{C})$ denotes a small finite group.

\end{theorem}

All the results that are well known in residue theory consider the hypothesis that the component of singular set is nondegenerate, see for instance \cite{Article-2}. The paper \cite{Article-17} provides a slight improvement of the results given by Baum-Bott by considering the degenerate case with restrictions. Let us consider $v$ a holomorphic vector field on $V$ and $W$ a component of singular set of $v$ such that the vector field is degenerate along of $W$. M. Dia proves the result of residue of $v$ at $W$ subject to the condition that there is a biholomorphism between a neighborhood of $W$ and a neighborhood of the zero section of the normal bundle of $W$. See Th\'eor\`eme A, B, C and D in \cite{Article-17}.

\subsection{Residues of logarithmic foliations of dimension one}

This section is dedicated to present some results about residues associated with logarithmic  holomorphic foliations of dimension one. For this we refer the readers to \cite{Article-0, Article-12.1,Article-14.3,  Article-20.11, Article-20.2} and the references therein.

The general index of a vector field tangent to hypersurfaces was defined and studied in terms of the homology of the complex of differential forms by X. Gomez-Mont, L. Giraldo and P. Mardesi\'c, see \cite{Article-20.11, Article-20.2}. The first result of this section deals with logarithmic index and was done by A. G. Aleksandrov in \cite{Article-0}. The author defines a logarithmic index using differential forms with logarithmic poles and proves that the homological index can be expressed via the logarithmic index.

Let $M$ be a complex manifold of dimension $n$, and let $\Omega_{M}^{q} , q \geq 0$, and $Der(M)$ be the sheaves of germs of holomorphic $q$-forms and vector fields on $M$, respectively. Let $D\subset M$ be a divisor which all of whose
irreducible components are of multiplicity one. Given $V \in Der(M)$ a vector field with an isolated singularity at a point $x \in D$ then the $\iota_{V}$-homology groups of the complex $(\Omega_{D,x}^{\bullet}, \iota_{V})$ are finite-dimensional vector spaces, where $\iota_{V} : \Omega_{M}^{q} \rightarrow \Omega_{M}^{q-1}$ is the interior multiplication (contraction). We can define the homological index of the vector field $V$ at the point $x \in D$ by

$$Ind_{hom,D,x}(V) := \sum_{i=0}^{n}(-1)^{i}\dim H_{i}(\Omega_{D,x},\iota_{V}).$$
Given a divisor $D$ we can consider the coherent analytic sheaves $\Omega_{M}^{q}(\log D)$, $q >0$  and $  Der_{M}(\log D) = T_{M}(-\log D)$ as in \cite{Article-0}. Take a vector field $V \in Der_{M} (\log D$). As above, the interior multiplication $\iota_{V}$ defines the structure of a complex on $\Omega_{M}^{\bullet}(\log D).$

\begin{lemma}\label{lemma: important}(\cite{Article-0}, Lemma 1, p.247) If all singularities of $V$ are isolated, then the $\iota_{V}$-homology groups of the complex $\Omega_{M}^{\bullet}(\log D)$ are finite-dimensional vector spaces.
\end{lemma}
It follows from Lemma \ref{lemma: important} that the \textit{logarithmic index} of the field $V$ is well defined at the point $x$,

$$ Ind_{\log D,x}(V) := \sum_{i=0}^{n}(-1)^{i}\dim H_{i}(\Omega_{D,x}(\log D),\iota_{V}).$$

These index are related bellow.

\begin{proposition}(\cite{Article-0}, Proposition 1, p. 248) Suppose that a point $x \in D$ is an isolated singularity of a vector field $V \in
Der(\log D)$, the germs $V_{i} \in  \mathcal{O}_{M,x}$ are determined by the expansion $V = \sum_{i} V_{i} \dfrac{\partial}{\partial z_{i}}$, and $J_{x}V = (V_{0},\ldots, V_{n})\mathcal{O}_{M,x}$. Then

$$Ind_{hom,D,x}(V) = \dim \mathcal{O}_{M,x}/J_{x}V - Ind_{\log D,x}(V).$$

\end{proposition}

In \cite{Article-14.3} the authors consider logarithmic foliation along $D$ and prove the residue formulas, namely, Baum-Bott type formulas for non-compact complex manifold, still considering the logarithmic vector field.

\begin{theorem}(\cite{Article-14.3}, Theorem 1, p. 6403) Let $\tilde{X}$ be an $n$-dimensional complex manifold such that $\tilde{X} = X-D$, where $X$ is an $n$-dimensional complex compact manifold and $D$ is a smooth
hypersurface on on $X$. Let $\mathcal{F}$ be a foliation of dimension one on $X$ with isolated singularities and logarithmic along $D$. Suppose that $Ind_{\log D,p}(\mathcal{F}) = 0$ for all $p \in \Sing(\mathcal{F}) \cap D$. Then

$$\int_{X} c_{n}(T_{X}(-\log D)-T_{\mathcal{F}}) = \sum_{p \in \Sing(\mathcal{F})\cap (X \setminus D)} \mu_{p}(\mathcal{F}).$$

\end{theorem}

In the same work, the authors consider that the divisor $D$ is a normal crossing hypersurface and one has

%In \cite{Article-14.3} Corr\^ea and Machado extended the Baum-Bott Theorem from compact complex manifolds to noncompact complex manifolds where the compactification is by an analytic divisor $D$ invariant by an one-dimensional holomorphic
%foliation $\mathcal{F}$ which is called by logarithmic foliation along $D$, in this works, the authors consider two cases, the first ones $D$ is a smooth
%hypersurface and the second ones $D$ a normally crossing hypersurface, see for instance the second ones

\begin{theorem}(\cite{Article-14.3}, Theorem 2, p. 6404) Let $\tilde{X}$ be an $n$-dimensional complex manifold such that $\tilde{X} = X-D$, where $X$ is an $n$-dimensional complex compact manifold, $D$ is a normally crossing hypersurface on $X$. Let $\mathcal{F}$ be a foliation on $X$ of dimension one, with isolated singularities (non-degenerates) and logarithmic along $D$. Then,

$$\int_{X} c_{n}(T_{X}(-\log D)-T_{\mathcal{F}}) = \sum_{p \in \Sing(\mathcal{F})\cap (\tilde{X})} \mu_{p}(\mathcal{F}).$$

\end{theorem}

In \cite{Article-12.1} the authors prove new versions of Gauss-Bonnet and Poincar\'e-Hopf theorems for complex $\partial$-manifolds of the type $\tilde{X}= X-D$, where $\dim X = n \geq 3$ and $D$ is a reduced divisor. More precisely,

\begin{theorem}(\cite{Article-12.1}, Theorem 1.1 p. 495)
Let $\tilde{X}$ be a complex manifold such that $\tilde{X} = X -  D$, where $X$ is an $n$-dimensional ($n\geq 3$)  complex compact manifold and $D$ is a reduced divisor on $X$. Given any (not necessarily irreducible) decomposition $D = D_1\cup D_2$, where $D_1$, $D_2$ have  isolated singularities and $C=D_1\cap D_2$  is a codimension $2$ variety and has isolated singularities,
\begin{itemize}
\item [ (i)] (Gauss-Bonnet type formula) the following  formula  holds\\ 

\small{
$$
\int_{X}c_{n}(\Omega^1_X(\log\,  D)) = (-1)^{n}\chi(\tilde{X})+\mu( D_1,S(D_1))+  \mu( D_2,S(D_2))-  \mu( C,S( C)).
$$}
\end{itemize}
\bigskip
\begin{itemize}
\item [ (ii)] (Poincar\'e-Hopf type formula) if $v$ is a holomorphic vector field on $X$, with isolated singularities  and logarithmic along $ D$,  we have that 

\small{
$$
 \chi(\tilde{X}) = PH(v, \Sing(v)) - GSV(v,  D_1,  S(v, D_1))  - GSV(v,  D_2,  S(v, D_2))  +  
 $$
 $$
+GSV(v,  C,  S(v, C))+(-1)^{n-1}\left[  \mu( D_1,S(D_1))+  \mu( D_2,S(D_2))-  \mu( C, S( C)) \right].
$$}

%$$\begin{array}{clll}
%\chi(\tilde{X}) &=& PH(v, \Sing(v)) - GSV(v,  D_1,  S(v, D_1))  - GSV(v,  D_2,  S(v, D_2)) + \\

%&  &  GSV(v,  C,  S(v, C))  + (-1)^{n-1}\left[  \mu( D_1,S(D_1))+  \mu( D_2,S(D_2))-  \mu( C, S( C)) \right].

%\end{array}$$

\end{itemize}
\end{theorem}

\subsection{Residues to flags}

In the following, we review some results on residues of flag of foliations which have emerged in recent years. For the background of flags, we refer to \cite{Article-9, Article-12.22, Article-20, Article-24} and references therein.

The next result (\cite{Article-9}, Theorem 2) serves two purposes.  The first one is to show that if $\mathcal{F} = (\mathcal{F}_{1}, \mathcal{F}_{2})$ is a flag of holomorphic foliations on $M$ then we can define the residue of $\mathcal{F}$ as an element of homology group. The second one is to show a Baum-Bott type theorem of the flag $\mathcal{F}$.

\begin{theorem}\label{1.0.14} Let $\mathcal{F} = (\mathcal{F}_{1}, \mathcal{F}_{2} )$ be a 2-flag of holomorphic foliations on a compact complex manifold
$M$ of dimension $n$. Let $\varphi_{1}, \varphi_{2}$ be homogeneous symmetric polynomials, respectively of degrees $d_{1}
$ and $d_{2}$, satisfying the Bott vanishing theorem to Flags. Then for each compact connected component $S$ of $\Sing(\mathcal{F})$ there exists a class, $\Res_{\varphi_{1}, \varphi_{2}} (\mathcal{F}, \mathcal{N}_{\mathcal{F}}; S ) \in H_{2n - 2(d_{1} + d_{2} )} (S; \mathbb{C})$, that we will call it of Baum-Bott Reisdue of Flag, such that
\begin{equation}\label{eq.0.2}
\sum_{\lambda}(\iota_{\lambda})_{\ast} \mbox{\Res}_{\varphi_{1}, \varphi_{2}} (\mathcal{F} , \mathcal{N}_{\mathcal{F}}; S_{\lambda} )=
(\varphi_{1}(\mathcal{N}_{12}). \varphi_{2}(\mathcal{N}_{2})) \frown [M]  
\end{equation}

\noindent in $H_{2n - 2(d_{1} + d_{2} )} (M; \mathbb{C})$, where $\iota_{\lambda}$ denotes the embedding of $S_{\lambda}$ in $M$.
\end{theorem}

	It is very interesting that we have an effective way to compute residues. The remainder of this section is devoted to this topic. First, we point out how the residues of two foliations in a flag are related (see \cite{Article-9}, Proposition 3, p. 1169).

\begin{theorem} For a flag $\mathcal{F} = (\mathcal{F}_{1}, \mathcal{F}_{2} )$ on $M$ with $\dim(\mathcal{F}_{1}) = \codim(\mathcal{F}_{2}) = 1$ and $\Sing(\mathcal{F}_{1}) \cap \Sing(\mathcal{F}_{2})$ admitting isolated singularities (only) we have
$$ \mbox{\Res}_{c_{n}} (\mathcal{F}_{2}, \mathcal{N}_{2}; p) = (-1)^{n} (n-1)! \mbox{\Res}_{c_{n}} (\mathcal{F}_{1}, \mathcal{N}_{1}; p), $$

\noindent where the residues involved are of the foliations $\mathcal{F}_{1}$ and $\mathcal{F}_{2}$.

\end{theorem}

 Let us introduce some notations. Given a flag $\mathcal{F} = (\mathcal{F}_{1}, \mathcal{F}_{2})$ on $M$ with  $\codim (\mathcal{F}_{i}) = k_{i}$ and $i =1,2,$ we denote by $\mbox{\Sing}_{k_{i}+1} \mathcal{F}_{i}$ the union of irreducible components of $\Sing(\mathcal{F}_{i})$ of pure codimension $k_{i}+1$. Take an irreducible component $Z \subset \mbox{\Sing}_{k_{1} + 1}(\mathcal{F}_{1})$ and a generic point $p \in Z$. Denote by $B_{p}$ a small ball centered at $p$ such that $S(B_{p}) \subset B_{p}$ is a sub-ball of dimension $n - k_{1} - 1$ (i.e., the same dimension of $Z$). 

Under the above setting (see \cite{book-0, book-4}) the de Rham class can be integrated over an oriented $(2k_{1} + 1)$-sphere $L_{p} \subset B_{p}^{\ast}$ and it defines the Baum-Bott residue of $\mathcal{F}$ at $Z$

$$ BB^{j}(\mathcal{F}; Z) := (2 \pi i)^{-k_{1} - 1} \int_{L_{p}} \theta^{12}\wedge (d \theta^{2})^{j} \wedge (d \theta^{12})^{k_{1} - j}\ \ \ \mbox{for each} \ \ \ 0 \leq j \leq k_{2}. $$

In particular, for the case where $\codim(\Sing(\mathcal{F})) \geq k_{1} + 1$, we get the formula to residue and the Baum-Bott theorem, see (\cite{Article-9}, Theorem 4, p.1173).

\begin{theorem} Let $\mathcal{F} = (\mathcal{F}_{1}, \mathcal{F}_{2})$ be a 2-flag of codimension $(k_{1}, k_{2})$ on a compact complex manifold $M$. If $\codim(\Sing(\mathcal{F})) \geq k_{1} + 1$, then for each $0 \leq j \leq k_{2}$ we have

$$c_{1}^{k_{1} - j + 1}(\mathcal{N}_{12}) \smile c_{1}^{j}(\mathcal{N}_{2}) = \sum_{Z \ \in \ \mbox{\Sing}_{k_{1}+1}(\mathcal{F}_{1}) \cup \mbox{\Sing}_{k_{1}+1}(\mathcal{F}_{2}) } \lambda_{Z}^{j}(\mathcal{F})[Z], $$

\noindent where $\lambda_{Z}^{j}(\mathcal{F}) = BB^{j}(\mathcal{F}, Z)$.

\end{theorem}

In 2020 Ferreira and Louren\c co in \cite{Article-20} extended the residue theory to flags of holomorphic distributions. The next result is about isolated singularities.

\begin{theorem}(\cite{Article-20}, Theorem 1.2)\label{2.11} Let $\mathcal{F} = (\mathcal{F}_{1}, \mathcal{F}_{2} )$ be a $2$-flag of holomorphic distributions on a compact complex manifold $M$ of dimension $n$, $\varphi_{1} \ and \ \varphi_{2}$ be homogeneous symmetric polynomials, respectively of degrees $d_{1}>0$ and $d_{2}> 0$ and $p$ an isolated point of $\Sing(\mathcal{F})$. Then

$$\Res_{\varphi_{1}, \varphi_{2}} (\mathcal{F} , \mathcal{N}_{\mathcal{F}}; p ) = 0.$$

\end{theorem}

 In the particular case where M is the projective space $\mathbb{P}^3$
	we obtain some advances in goal of compute residue of flags.
	Furthermore, if the singular scheme of flag has only one irreducible component the following result gives us an effectively way to compute residues, see ([36], Theorem 1.3).

\begin{theorem}\label{1.2} Let $\mathcal{F} = (\mathcal{F}_{1}, \mathcal{F}_{2} )$ be a $2$-flag of holomorphic foliations on $\mathbb{P}^{3}$ with $\deg(\mathcal{F}_{i}) = d_{i}$ thus

\begin{equation}\label{2}
(1 + d_{1} - d_{2})\sum_{Z \in S_{1}(\mathcal{F}_{2})}\deg(Z)\Res_{\varphi_{2}}(\mathcal{F}_{2}|_{B_{p}}; p)=\sum_{Z \in S_{1}(\mathcal{F})} \Res_{c_{1}\varphi_{2}}(\mathcal{F}, \mathcal{N}_{\mathcal{F}}; Z),
\end{equation}

\noindent where $\deg(Z)$ is the degree of the irreducible component $Z$, $\Res_{\varphi_{2}}(\mathcal{F}_{2}|_{B_{p}}; p)$ represents the Grothendieck residue of the foliation $\mathcal{F}_{2}|_{B_{p}}$ at $\{p\} = Z \cap B_{p}$  with $B_{p}$ a transversal ball and either $\varphi_{2}=c_{1}^{2}$ or $\varphi_{2}=c_{2}.$

\end{theorem} 
We believe that if we work with the theory of residue currents, see \cite{Article-25}, we can present a more general expression to compute residue, mainly for the degenerate case.

\section{Some applications}

\subsection{Residues and the Poincar\'e problem}

Motivated by Darboux's works and the problem of the algebraic integrability of differential equations, in 1891 H. Poincar\'e formulated the following question, see \cite{Article-24.1}
\\

\noindent \textit{"Is it possible to bound the degree of an irreducible curve such that is invariant by a foliation in terms of the degree of foliation?"} \\

\noindent This problem is equivalent to  asking when a holomorphic foliation on $\mathbb{P}^{2}$ admits a rational first integral. This question is known as the \textit{Poincar\'e problem}. It is well known that such a bound does not exist in general, see \cite{Article-11.1}. However, under certain hypotheses, there are several partial answers and generalizations even for flags and for Pfaff systems; see for instance \cite{Article-9.1, Article-9.2, Article-11.1,  Article-12.2, Article-12.22,Article-14.2, Article-17.1,Article-20,  Article-20.12, Article-25.1, Article-27.1}. The residue theory, and especially the Baum-Bott Theorem, are powerful tools and obstructions for several problems related to foliations with singularities.

One hundred years after, in (\cite{Article-11.1}, Theorem 1, p.891) Cerveau and Lins Neto gave a first partial answer to Poincar\'e problem.  For this, consider 
$S$ a projective nodal curve, that is, all its singularities are of normal crossing type, with reduced homogeneous equation $f=0$, of degree $m$.

\begin{theorem} Let $\mathcal{F}$ be a foliation in $\mathbb{CP}(2)$ of degree $n$, having
$S$ as separatrix. Then $m \leq n+2$. Moreover if $m=n+2$ then $f$ is reducible and $\mathcal{F}$ is of logarithmic type, that is given by a rational closed
form $\sum \lambda_{i} \dfrac{df_{i}}{f_{i}}$ where $\lambda_{i} \in \mathbb{C}$ and the $f_{i}$ are homogeneous polynomials.

\end{theorem}

Some years later, Carnicer in (\cite{Article-9.2}, Theorem, p.289) presents the same quota to the Poincar\'e problem with other hypothesis, namely, the foliation does not have dicritical singularities into the curve. Roughly speaking, a  foliation $\mathcal{F}$ is dicritical if, and only if, there exists an irreducible surface $Z$ such that restriction of $\mathcal{F}$ to $Z$ contains infinitely many distinct separatrices.

\begin{theorem} Let $\mathcal{F}$ be a foliation of $\mathbb{P}^{2}$ and let $C$ be an algebraic curve in
$\mathbb{P}^{2}$. Suppose that $C$ is invariant by $\mathcal{F}$ and there are no dicritical singularities
of $\mathcal{F}$ in $C$. Then

$$\deg C \leq \deg \mathcal{F} + 2.$$

\end{theorem}

M. Soares in (\cite{Article-27.1}, Theorem B, p.496) extended the Poincar\'e problem from $\mathbb{P}^{2}$ to $\mathbb{P}^{n}$. Using the  Baum-Bott theorem, Soares improved the above quota. In fact, let $V \hookrightarrow  \mathbb{P}^{n}, n \geq 2,$ be an irreducible non-singular algebraic hypersurface of degree $d_{0}$ invariant by $\mathcal{F}^{d}$, a non-degenerated one-dimensional holomorphic foliation of degree $d \geq 2$. Then we have.

\begin{theorem} $$d_{0} \leq d + 1.$$

\end{theorem}

In \cite{Article-9.01} Brunella says that the GSV-index is the obstruction to a positive solution to Poincar\'e problem and gives a simple condition that implies the  non-negativity of this index. Motivated by Brunella's work, Corr\^ea and Machado in \cite{Article-14.2}  introduced a GSV type index for invariant varieties by holomorphic Pfaff systems on projective manifolds. The authors  proved, with certain hypotheses, a non-negativity property for this index. As a consequence we have the following result.

\begin{theorem} Let $\omega \in H^{0}(\mathbb{P}^{n}, \Omega^{k}_{\mathbb{P}^{n}}(d+k+1))$ be a holomorphic Pfaff system
of rank $k$ and degree $d$. Let $V \subset \mathbb{P}^{n}$ be a reduced complete intersection variety,
of codimension $k$ and multidegree $(d_{1},\ldots, d_{k})$, invariant by $\omega$. Suppose that $\Sing(\omega, V )$ has codimension one in $V$, then

$$\sum_{i}GSV(\omega,V, S_{i}) \deg(S_{i}) = [d+k+1-(d_{1}+ \dots +d_{k})](d_{1}\dots d_{k}),$$

\noindent where $S_{i}$ denotes an irreducible component of $\Sing(\omega, V)$. In particular, if $GSV(\omega, V, S_{i}) \geq 0$, for all $i$, we have

$$d_{1} +\dots + d_{k} \leq d + k + 1.$$

\end{theorem}

\subsection{Residues and the non-existence of minimal sets}

Let $\mathcal{F}$ be a holomorphic foliation on a compact manifold $X.$ A compact non-empty subset $\mathcal{M} \subset X$ is said to be a  \textit{nontrivial minimal set} for $\mathcal{F}$ if the following properties
are satisfied

\begin{enumerate}
	
	\item[a)] $\mathcal{M}$ is invariant by $\mathcal{F};$
	
	\item[b)] $\mathcal{M} \cap \Sing(\mathcal{F}) = \emptyset;$
	
	\item[c)] $\mathcal{M}$ is minimal with respect to properties $a)$ and $b)$.
	
\end{enumerate}

Next, we shall see how the residue theory can be used in the problem of existence or not of nontrivial minimal sets for foliations.

 The above problem for codimension one holomorphic foliations on complex projective spaces was considered by Camacho, Lins Neto and Sad in \cite{Article-9.22}. The authors addresses the problem on $\mathbb{P}^{2}$. More precisely, they prove a geometrical property of minimal sets, i.e., by applying the Maximum Principle for harmonic functions, they prove that $\mathcal{F}$ has at most one non-trivial minimal set. Moreover, under generic conditions imposed on the singularities of foliation, all leaves accumulate on that set. Anyway, in general, the question on $\mathbb{P}^{2}$ remains as an open problem.

  In the case where the complex projective space has higher dimensions, Lins Neto (\cite{Article-23.1}, Theorem 1, p. 1370) extended the above discussion and prove the following.

\begin{theorem}\label{alcideSO} Codimension 1 foliations on $\mathbb{P}^{n}, n\geq 3$, have no nontrivial minimal sets.
\end{theorem}

 We would like to emphasize that the residue theory is essential to prove Theorem \ref{alcideSO}. In fact, one of the main arguments in the proof is the following proposition (\cite{Article-23.1}, Proposition 1, p. 1372) which is a consequence of a Baum-Bott theorem. 

\begin{proposition}
Let $\mathcal{F}$ be a foliation of degree $k$ on $\mathbb{CP}^{2}$, with singular set $\Sing(\mathcal{F})$ of complex codimension $2$. Then
\[
\sum_{p \, \in\, \Sing(\mathcal{F})} BB(\mathcal{F},p) = (k+2)^{2}.
\] 
In particular $\displaystyle{\sum_{p \, \in\, \Sing(\mathcal{F})} BB(\mathcal{F},p)}$ is positive for any foliation $\mathcal{F}$ on $\mathbb{CP}^2$.
\end{proposition}

In (\cite{Article-9.11}, Theorem 1.2, p.296] the authors, by using Baum-Bott Theorem, have generalized Lins Neto's result for codimension one holomorphic foliations on complex projective manifolds and such that the Picard group is cyclic.

\begin{theorem}\label{burnella} Let $X$ be a complex projective manifold of dimension at least $3$ and with Pic$(X) = \mathbb{Z}$, and let $\mathcal{F}$ be a codimension one foliation
on $X$. Then every leaf $L$ of $\mathcal{F}$ accumulates to $\Sing(\mathcal{F})$:

$$\overline{L} \cap \Sing(\mathcal{F}) \neq \emptyset.$$

\end{theorem}

Theorem \ref{burnella} is a partial answer to Brunella's conjecture, see (\cite{Article-9.001}, Conjecture 1.1, p.3102), in the special case when Pic$(X) = \mathbb{Z}$. Recently, in 2021 M. Adachi and J. Brinkschulte have proved the Brunella's conjecture without hypotheses on manifold $X$, see (\cite{Article-1.1}, Main Theorem, p.1).

\begin{theorem} Let $X$ be a compact complex manifold of dimension $\geq 3$. Let $\mathcal{F}$ be a codimension one holomorphic foliation on $X$ with ample normal bundle $N_{\mathcal{F}}$. Then every leaf of $\mathcal{F}$ accumulates to $\Sing(\mathcal{F})$.
\end{theorem}

Brunella's conjecture has been stated for higher codimensional foliations (\cite{Article-14.22}, Conjecture 1.2, p. 1236) and still remains as an open problem. However, in \cite{Article-14.222} the authors by using a Brunella-Khanedani-Suwa variational type residue theorem for invariant currents by holomorphic foliations, provide conditions for the accumulation of the leaves to the intersection of the singular set of a holomorphic foliation with the support of an invariant current.

\end{document}